\documentclass[11pt]{article}
\usepackage{amsmath}
\usepackage{amsfonts}
\usepackage{amsthm}
\usepackage{pdfpages}
\def \dsp {\displaystyle}

\newtheorem{theo}{Theorem}[section]

\newcommand{\eProof}{\hfill { \mbox{$\Box$}}\medskip \smallskip} 

\def \dsp {\displaystyle}

\newcommand{\C}{\mathbb{C}}
\newcommand{\D}{\mathbb{D}}

\newcommand{\T}{\mathbb{T}}

\renewcommand{\Re}{\mathop{\rm Re}}
\renewcommand{\Im}{\mathop{\rm Im}}

\begin{document}
\setlength{\parskip}{0.0cm}

\title{Two families of orthogonal polynomials on the unit circle from basic hypergeometric functions\thanks{This  work has received  support from the funding bodies CNPq, Brazil (grant  305073/2014-1) and FAPESP, Brazil (grants 2016/09906-0 and 2017/12324-6).}}
\author{{A. Sri Ranga } \\[1.5ex] 
{\small DMAp, IBILCE, } 
{\small UNESP - Universidade Estadual Paulista} \\
{\small 15054-000  S\~{a}o Jos\'{e} do Rio Preto,  SP,  Brazil} \\[1ex]}
\date{\today}
\maketitle 


\begin{abstract} 
\noindent The sequence $\{\,_2\phi_1(q^{-k},q^{b+1};\,q^{-\overline{b}-k+1};\, q, q^{-\overline{b}+1/2} z)\}_{k \geq 0}$ of basic hypergeometric polynomials is known to be orthogonal on the unit circle with respect to the weight function $|(q^{1/2}e^{i\theta};\,q)_{\infty}/(q^{b+1/2}e^{i\theta};\,q)_{\infty}|^2$. This result, where one must take the parameters $q$ and $b$ to be $0 < q < 1$ and $\Re(b) > -1/2$, is due to P.I. Pastro \cite{Pastro-1985}.  In the present manuscript we deal with the  orthogonal polynomials $\hat{\Phi}_{n}(b;.)$ and $\check{\Phi}_{n}(b;.)$ on the unit circle with respect to the two parametric families of weight functions $\hat{\omega}(b; \theta) = |(e^{i\theta};\,q)_{\infty}/(q^{b}e^{i\theta};\,q)_{\infty}|^2$ and $\check{\omega}(b;\theta) = |(qe^{i\theta};\,q)_{\infty}/(q^{b}e^{i\theta};\,q)_{\infty}|^2$, where $0 < q < 1$ and $\Re(b) > 0$. With the use of the basic hypergeometric polynomials $ _2\phi_1(q^{-k},q^{b};\,q^{-\overline{b}-k+1};\, q, q^{-\overline{b}+1} z)$, $k \geq 0$, which have zeros on the unit circle when $\Re(b) > 0$, simple expressions for the (monic) polynomials $\hat{\Phi}_{n}(b;.)$ and $\check{\Phi}_{n}(b;.)$, their norms, the associated Verblunsky coefficients and also the respective Szeg\H{o} functions are found. 

\end{abstract}


\vspace{2ex}

{\noindent}Keywords: 
Orthogonal polynomials on the unit circle, 
Basic hypergeometric functions, Positive chain sequences.  \\

{\noindent}2010 Mathematics Subject Classification: 42C05, 33C45.

\setcounter{equation}{0}
\section{Introduction}

Orthogonal polynomials on the unit circle (in short, OPUC) are important objects of study in Classical Analysis. Like their counterpart on the real line, they have deep connections and applications in many areas of Mathematics and Engineering. Their systematic study, started by Szeg\H{o}  and Geronimus (see \cite{Szego-Book} and \cite{Geronimus-AMSTransl-1977}), still remains very active, especially because of their applications in spectral theory \cite{Simon-Book-p2, Simon-2011}.

Given a nontrivial probability measure $\mu$ on the unit circle $\T: = \{\zeta=e^{i\theta}\!\!: \, 0 \leq \theta \leq 2\pi \}$ the associated orthonormal OPUC $\varphi_k(z)=\kappa_k z^k + lower\ degree\ terms$, $k \geq 0$,  are defined by  $\kappa_k>0$ and 
\[
    \int_{\T}  \overline{\varphi_{j}(\zeta)}\, \varphi_{k}(\zeta)\, d \mu(\zeta) = \int_{0}^{2\pi} \overline{\varphi_{j}(e^{i\theta})}\, \varphi_{k}(e^{i\theta})\, d \mu(e^{i\theta}) = \delta_{j,k}, \quad j,k=0,1,2, \ldots, 
\]
where $\delta_{j,k}$ stands for the Kronecker delta. Among their fundamental properties is that all their zeros belong to the open unit disk $\D:= \{ z\in \C: \, |z|<1\} $, and that they satisfy the Szeg\H{o} recurrence,
\begin{equation} \label{Szego-A-RR}
   \Phi_{k}(z) =  z \Phi_{k-1}(z) - \overline{\alpha}_{k-1}\,      \Phi_{k-1}^{\ast}(z), \quad 
k \geq 1,
\end{equation}
given here in terms of the monic OPUC $\Phi_{k}(z) = \varphi_{k}(z)/\kappa_k$, $k \geq 0$. The coefficients  $\alpha_{k-1} = - \overline{\Phi_{k}(0)}$ are known as the Verblunsky coefficients, and $\Phi_{k}^{\ast}(z) = z^{k} \overline{\Phi_{k}(1/\bar{z})}$.  It is well known that  $|\alpha_k| < 1$ for $k \geq 0$, and that the sequence  $\{\alpha_k \}_{k\geq 0}$ uniquely determines the  measures $\mu$ on $\T$ (see e.g.~\cite{Simon-Book-p1}, as well as \cite{ENZG1}). 

If the theory of OPUC still leaves something to be highly desired, perhaps it is about having more concrete examples (examples with explicit formulas) of these polynomials. Most of the known examples with explicit formulas come from  orthogonal polynomials on the interval $[-1,1]$ via the Szeg\H{o}-transformation (see \cite{Szego-Book}) or via the DG-transformation (see \cite{Zhedanov-JAT1998}) and, hence, in these examples the values of associated Verblunsky coefficients are restricted to be real. Among the few known examples of monic OPUC with complex Verblunsky coefficients, the ones that we like to highlight are the following. 

With $\Re(b) > -1/2$, if $\Phi_k(z) = \frac{(b+\bar{b}+1)_k}{(b+1)_k}\ _2F_1(-k,b+1;b+\bar{b}+1;1-z)$, $k \geq 0$, then  $\{\Phi_k\}_{k \geq 0}$ is the sequence of monic OPUC such that    
\[
   \begin{array}l
   \dsp\frac{2^{b+\bar{b}}\,|\Gamma(b+1)|^2}{2\pi\, \Gamma(b+\bar{b}+1)}
        \int_{0}^{2\pi} \overline{\Phi_{j}(e^{i\theta})}\, \Phi_{k}(e^{i\theta})\, (e^{\pi-\theta})^{\Im(b)} (\sin^{2}(\theta/2))^{\Re(b)} d\theta \\[1ex]
       \dsp \hspace{53ex}   = \frac{(b+\bar{b}+1)_k\,k!}{ |(b+1)_k|^2}\, \delta_{j,k}.
   \end{array}
\]
We have for the associated Verblunsky coefficients  $\alpha_{k-1}  = -(b)_k/ (\overline{b}+1)_k$, $ k \geq 1$.    
For information regarding the definitions and properties of the Pochhammer symbols $(.)_{k}$ and the hypergeometric functions $_2F_1$,  we refer to \cite{AndAskRoy-book}. 

The above parametric family of OPUC came to be of public knowledge in \cite{Ranga-PAMS2010}.  However, it is important to mention that this family of OPUC are a subfamily of a family of biorthogonal polynomials presented by Askey in the Gabor Szeg\H{o}: Collected papers \cite[p.\,304]{Askey-SzegoCollected1-1982}. 

The second parametric family of OPUC with complex Verblunsky coefficients, attributed to Pastro \cite{Pastro-1985},  that we like to mention here is   
\begin{equation} \label{Eq-pastro-Polys}
  \begin{array}{ll}
     \Phi_{k}^{(b)}(z) = \frac{(q^{\overline{b}};\,q)_k}{(q^{b+1};\,q)_k}\, q^{k/2}
                          \,_2\phi_1\Big(\begin{array}{c}
                             q^{-k},\, q^{b+1} \\
                              q^{-\overline{b}-k+1}
                             \end{array}\!\!; \, q,\, q^{-\overline{b}+1/2}z\Big), \quad k \geq 0,
  \end{array}
\end{equation}
again with  $\Re(b) > -1/2$.  These monic polynomials satisfy the orthogonality
\[
     \int_{\mathbb{T}} \overline{\Phi_{j}^{(b)}(\zeta)}\, \Phi_{k}^{(b)}(\zeta)\,  d \mu^{(b)}(\zeta) =  \rho_{k}^{(b, b+\overline{b}-1)}  = \frac{(q;\,q)_k\, (q^{b+\overline{b}+1};\,q)_k}{(q^{b+1};\,q)_k\,  (q^{\overline{b}+1};\,q)_k}\,\delta_{j,k}, \quad
\]
where the probability measure $\mu^{(b)}$ on the unit circle is such that 
\[
     d \mu^{(b)}(\zeta) = \frac{(q;\,q)_{\infty}(q^{b+\overline{b}+1};\,q)_{\infty}} {(q^{b+1};\,q)_{\infty}(q^{\overline{b}+1};\,q)_{\infty}}\, \frac{|(q^{1/2}\zeta;\,q)_{\infty}|^2} {|(q^{b+1/2}\zeta;\,q)_{\infty}|^2} \frac{1}{2\pi i\zeta}\, d\zeta .
\]

For information on the definitions and properties of  q-Pochhammer symbols $(\,.\,;q)_{k}$ and  basic hypergeometric (or q-hypergeometric) functions $_2\phi_1$ we refer, for example, to Gasper and Rahman \cite{GasRah-book} and Koekoek and Swarttouw \cite{KoeSwa-book}.  

Apart from the above two examples of OPUC with complex Verblunsky coefficients, it is important that we mention the system of OPUC with constant Verblunsky coefficients. This system of OPUC, also known as Geronimus polynomials,  has been thoroughly studied by many including Geronimus \cite{Geronimus-AMSTransl-1977} and  Golinskii, Nevai and  Van Assche \cite{GolinNevaiAssche-JAT1995}  (see also  Simon \cite[p.\,83]{Simon-Book-p1}).  

Our present aim is to consider the OPUC and related information  associated with the positive measures $\hat{\mu}^{(b)}$ and $\check{\mu}^{(b)}$ on the unit circle given by 
\begin{equation} \label{Eq-NewMeasure-1}
     d \hat{\mu}^{(b)}(\zeta) =  \hat{\sigma}^{(b)}\, \frac{|(\zeta;\,q)_{\infty}|^2} {|(q^{b}\zeta;\,q)_{\infty}|^2}\frac{1}{2 \pi i \, \zeta} \,d\zeta  
\end{equation}
and
\begin{equation} \label{Eq-NewMeasure-2}
    d \check{\mu}^{(b)}(\zeta) =  \check{\sigma}^{(b)}\, \frac{|(q\zeta;\,q)_{\infty}|^2} {|(q^{b}\zeta;\,q)_{\infty}|^2}\frac{1}{2 \pi i \, \zeta} \,d\zeta, 
\end{equation}
where $\mathcal{R}e(b) > 0$. The values of the positive constants $\hat{\sigma}^{(b)}$ and $\check{\sigma}^{(b)}$, so that the measures $\hat{\mu}^{(b)}$ and $\check{\mu}^{(b)}$ are probability measures, are given respectively in  \eqref{Eq-tau1} and \eqref{Eq-tau2}. 

The objects that play important roles in this manuscript are the modified basic hypergeometric polynomials 
\begin{equation} \label{Eq-Rk-ExplicitForm}
   R_{k}(b;\, z)   =  \frac{(q^{\overline{b}};\,q)_k}{(q^{\lambda}\cos(\eta_q);\,q)_k}\,
   \,_2\phi_1\Big(\begin{array}{c}
                   q^{-k},\, q^{b} \\
                   q^{-\overline{b}-k+1}
                  \end{array}\!\!;\, q,\, q^{-\overline{b}+1}z\Big), \quad  k \geq 0,
\end{equation}
and the sequence  $\{d_{k+1}^{(b)}\}_{k \geq 1}$, where 
\begin{equation} \label{TTRR-q-Rn-dn-coeffs}
      d_{k+1}^{(b)} = \frac{(1-q^{k})\, (1-q^{2 \lambda+k-1})}{4(1 - q^{\lambda+k-1}\cos(\eta_q))\,(1 - q^{\lambda+k}\cos(\eta_q))}, \quad k \geq 1. 
\end{equation}
Here,  $b = \lambda - i \eta$ and $\eta_{q} = \eta \ln(q)$.  As detailed in Section \ref{Sec-SelfInversive-polynomials} of this manuscript, with the choice $Re(b) = \lambda > 0$ the sequence  $\{d_{k+1}^{(b)}\}_{k \geq 1}$ is a positive chain sequence.  For a good source of information on positive chain sequences we cite Chihara \cite{Chihara-book}. 

The  multiplication (or modification) factor $(q^{\overline{b}};\,q)_k/(q^{\lambda}\cos(\eta_q);\,q)_k)$ in \eqref{Eq-Rk-ExplicitForm} is  such that there hold $R_{k}(b;\, z) = R_{k}^{\ast}(b;\, z)$, $k \geq 0$  and, further,  the sequence $\{R_{k}(b;\, .)\}_{k\geq 0}$ satisfies the  nice three term recurrence formula  
\begin{equation} \label{TTRR-q-Rn}
   R_{k+1}(b;\, z) = \big[(1 + i\,c_{k+1}^{(b)})z + (1 - i\,c_{k+1}^{(b)})\big]R_{k}(b;\, z) - 4 d_{k+1}^{(b)} z R_{k-1}(b;\, z), 
\end{equation}
for $k \geq 1$, with $R_{0}(b;\, z) = 1$ and $R_{1}(b;\, z) = (1 + i\,c_{1}^{(b)})z + (1 - i\,c_{1}^{(b)})$, where
\begin{equation} \label{TTRR-q-Rn-coeffs}
   c_{k}^{(b)} = \frac{q^{\lambda+k-1}\sin(\eta_q)}{1 - q^{\lambda+k-1}\cos(\eta_q)}, \quad k \geq 1.
\end{equation}

Our main result  with respect to the measure $\hat{\mu}^{(b)}$ is the following theorem, the proof of which is given in Section  \ref{Sec-New-OPUC1}. 

\begin{theo} \label{Thm-New-OPUC1}
Let $b = \lambda-i\eta$, $\eta_q = \eta\, \ln(q)$ and $\lambda > 0$. Then the sequence $\{\hat{\Phi}_{k}(b; z)\}_{k \geq 0}$ of monic OPUC  with respect to the positive measure  $\hat{\mu}^{(b)}$ given by \eqref{Eq-NewMeasure-1} is such that   
\[
    \hat{\Phi}_{k}(b; z) = \frac{(q^{\lambda}\cos(\eta_{q}); q)_{k+1}}{(q^{b};q)_{k+1}}\, \frac{R_{k+1}(b;z) - 2(1-\ell_{k+1}^{(b)}) R_{k}(b;z)}{z-1}, \quad k \geq 0.
\]
Here, $\ell_{k+1}^{(b)}$, $k \geq 0$ are such that $\ell_{1}^{(b)} = 0$ and $\ell_{k+1}^{(b)} = d_{k+1}^{(b)}/ (1-\ell_{k}^{(b)})$, $k \geq 1$.
In particular, the associated Verblunsky coefficients satisfy  
\[
    \hat{\alpha}_{k-1}^{(b)} = -   \Big[ 1 - 2 \ell_{k+1}^{(b)} \frac{1- q^{\lambda+k} \cos(\eta_{q})}{1 - q^{\overline{b}+k}}\Big] \frac{(q^{b};q)_{k}}{(q^{\overline{b}};q)_{k}}, \quad k \geq 1.
\]
Moreover, if  $\hat{\mu}^{(b)}$ is a probability measure  and  if $\hat{\phi}_{k}(b;z) = \hat{\kappa}_k^{(b)}\, \hat{\Phi}_{k}(b; z)$ are the associated orthonormal polynomials then
\[
     [\hat{\kappa}_{k}^{(b)}]^{-2} = \frac{(q;q)_{k}\,(q^{2\lambda};q)_{k}}{(q^{b+1};q)_{k}\,(q^{\overline{b}+1};q)_{k}} \frac{1-q^{\lambda+k}\cos(\eta_{q})}{1-q^{\lambda}\cos(\eta_{q})} (1 - \ell_{k+1}^{(b)}), \quad  k \geq 0.
\] 
\end{theo}

The sequence $\{\ell_{k+1}^{(b)}\}_{k \geq 0}$ is  the minimal parameter sequence of the positive chain sequence $\{d_{k+1}^{(b)}\}_{k \geq 1}$. An explicit expression for $\ell_{k+1}^{(b)}$ for any $k \geq 1$ is also given in Section \ref{Sec-New-OPUC1}. Finally, in Section \ref{Sec-New-OPUC1} the Szeg\H{o} function associated with the measure $\hat{\mu}^{(b)}$ is also explicitly found.

Now with respect to the measure $\check{\mu}^{(b)}$ our main result  is the following theorem.   

\begin{theo} \label{Thm-New-OPUC2}
Let $b = \lambda-i\eta$, $\eta_q = \eta\, \ln(q)$ and $\lambda > 0$. Then the sequence $\{\check{\Phi}_{k}(b; z)\}_{k \geq 0}$ of monic OPUC  with respect to the positive measure  $\check{\mu}^{(b)}$ given by \eqref{Eq-NewMeasure-2} is such that   
\[
    \check{\Phi}_{k}(b; z) = \frac{(q^{\lambda}\cos(\eta_{q}); q)_{k}}{(q^{b};q)_{k}}\, [R_{k}(b;z) - 2(1-M_k^{(b)})R_{k-1}(b;z)], \quad k \geq 1.
\]
Here, $M_{k}^{(b)}$, $k \geq 1$, are such that  $M_{k+1}^{(b)} = d_{k+1}^{(b)}/ (1-M_{k}^{(b)})$, $k \geq 1$, with 
\[ 
    M_{1}^{(b)} = \frac{1}{2}\frac{1 - q^{b}}{1 - q^{\lambda}\cos(\eta_{q})}\Big(1 - \frac{\int_{\T}\zeta\, d \check{\mu}^{(b)}(\zeta)}{\int_{\T} d \check{\mu}^{(b)}(\zeta)}\Big).
\]
The associated Verblunsky coefficients satisfy  
\[
    \check{\alpha}_{k-1}^{(b)} = \Big[ 1 - 2 M_{k}^{(b)} \frac{1- q^{\lambda+k-1} \cos(\eta_{q})}{1 - q^{\overline{b}+k-1}}\Big] \frac{(q^{b};q)_{k-1}}{(q^{\overline{b}};q)_{k-1}}, \quad k \geq 1.
\]
Moreover, if $\check{\mu}^{(b)}$ is a probability measure and  if $\check{\phi}_{k}(b;z) = \check{\kappa}_k^{(b)}\, \check{\Phi}_{k}(b; z)$ are the associated orthonormal polynomials then
\[
     [\check{\kappa}_{k}^{(b)}]^{-2} = \frac{(q;q)_{k}\,(q^{2\lambda};q)_{k}}{(q^{b};q)_{k}\,(q^{\overline{b}};q)_{k}}  \frac{1-q^{\lambda}\cos(\eta_{q})}{1-q^{\lambda+k}\cos(\eta_{q})}\frac{M_1^{(b)}}{M_{k+1}^{(b)}},\quad k \geq 0.
\] 
\end{theo}

The proof of Theorem \ref{Thm-New-OPUC2} is given in Section  \ref{Sec-New-OPUC2} of this manuscript.
As shown also  in Section \ref{Sec-New-OPUC2},  the sequence $\{M_{k+1}^{(b)}\}_{k \geq 0}$ is the maximal parameter sequence of the positive chain sequence $\{d_{k+1}^{(b)}\}_{k \geq 1}$. Explicit expression for $M_{k+1}^{(b)}$ for any $k \geq 0$ is also given in this section.  Finally in this section, the Szeg\H{o} function associated with the measure $\check{\mu}^{(b)}$ is also explicitly found.

The manuscript is organized as follows.  In Section  \ref{Sec-SelfInversive-polynomials} we provide some results  connected with the family $R_k(b;.)$ of basic hypergeometric  polynomials.  These results form the building blocks for the results obtained in Sections \ref{Sec-New-OPUC1} and \ref{Sec-New-OPUC2}. Specifically, Section \ref{Sec-New-OPUC1} deals with  the OPUC with respect to the measure $\hat{\mu}^{(b)}$ and Section \ref{Sec-New-OPUC2} gives information about the OPUC with respect to the measure $\check{\mu}^{(b)}$. 

\setcounter{equation}{0}
\section{Polynomials with zeros on the unit circle} \label{Sec-SelfInversive-polynomials}

For $0 < q < 1$ and $b \neq 0, -1, -2, \ldots$, we now consider the family of monic   polynomials given by  
\begin{equation} \label{Eq-P_k(b)}
   P_{k}(b;\, z) = \frac{(q^{\overline{b}};\,q)_k}{(q^{b};\,q)_k}\,
   \,_2\phi_1\Big(\begin{array}{c}
                   q^{-k},\, q^{b} \\
                   q^{-\overline{b}-k+1}
                  \end{array}\!\!;\, q,\, q^{-\overline{b}+1}z\Big), \quad k \geq 0.
\end{equation}
This is the subfamily $B_{k}^{(b-1, b+\overline{b}-2,b)}(z)$ of the family of basic hypergeometric polynomials 
\[
     B_{k}^{(b, c,d)}(z) = \frac{(q^{c-b+1};\,q)_k}{(q^{b+1};\,q)_k}\, q^{k(b-d+1)}
                  \,_2\phi_1\Big(\begin{array}{c}
                             q^{-k},\, q^{b+1} \\[0.5ex]
                              q^{-c+b-k}
                             \end{array}\!\!;  \, q,\, q^{-c+d-1}z\Big), \quad k \geq 1,
\]
studied in \cite{CosGodLamRan-2012}. Hence, from results obtained in \cite{CosGodLamRan-2012},  
\begin{equation} \label{TTRR-q-Pn}
   \begin{array} l
     P_{k+1}(b;\, z) = (z + \mathfrak{C}_{k+1}^{(b)})P_{k}(b;\, z) - \mathfrak{D}_{k+1}^{(b)}\,z\, P_{k-1}(b;\, z), 
   \end{array}
\end{equation}
for $k \geq 1$, with $P_{0}(b;\, z) = 1$ and $P_{1}(b;\, z) = z + \mathfrak{C}_{1}^{(b)}$, where
\begin{equation*} \label{TTRR-Coeff-q-Rn-Temp}
    \mathfrak{C}_k^{(b)}  = \frac{1-q^{\overline{b}+k-1}}{1-q^{b+k-1}}, \quad
    \mathfrak{D}_{k+1}^{(b)} =  \frac{(1-q^{k})\, (1-q^{b+\overline{b}+k-1})}{(1-q^{b+k-1})\, (1-q^{b+k})},\quad k \geq 1.
\end{equation*}
Moreover, there holds  the orthogonality
\begin{equation} \label{Orthogonality-general-Lb}
     \mathcal{L}^{(b)}[\zeta^{-j}P_{k}(b;\, \zeta)] = \delta_{k,j}\, \rho_k^{(b)}, \quad 0 \leq j \leq k, \quad k \geq 1,
\end{equation}
with respect to the quasi-definite moment functional
\begin{equation} \label{Explict-Exp-Moments-q-Rn}
\mathcal{L}^{(b)}[\zeta^{-j}] = \frac{(q^{-b+1};\,q)_j}{(q^{\overline{b}+1};\,q)_j}\, q^{jb}, \quad j = 0, \pm 1, \pm 2, \ldots \ .
\end{equation}
Here, \ $\rho_k^{(b)} = \frac{(q;\,q)_k\, (q^{b+\overline{b}};\,q)_k}{(q^{b};\,q)_k\,  (q^{\overline{b}+1};\,q)_k}$.

Let $b = \lambda - i \eta$ and $\eta_{q} = \eta\, \ln(q)$.  Then with the observation that $1 - q^{b+k} = 1 - q^{\lambda+k} \cos(\eta_q) + i\, q^{\lambda+k} \sin(\eta_q)$, let us consider the sequence of polynomials $\{R_{k}(b;\, .)\}_{k \geq 0}$ given by $R_{k}(b;\, z) = \frac{(q^{b};\,q)_k}{(q^{\lambda}\cos(\eta_q);\,q)_k}\, P_{k}(b;\, z)$, $k \geq 0$.

 Then the polynomials $R_{k}(b;\, z)$ takes the form \eqref{Eq-Rk-ExplicitForm}. 
From (\ref{TTRR-q-Pn}) one can easily verify that these polynomials also satisfy the three term recurrence formula \eqref{TTRR-q-Rn}. 

 From \eqref{TTRR-q-Rn} and \eqref{TTRR-q-Rn-coeffs} 
\[
    R_{k}(b;\, 0) = \prod_{j=1}^{k} (1-ic_{j}^{(b)}) = \frac{(q^{\overline{b}}; q)_{k}}{(q^{\lambda}\cos(\eta_{q}); q)_{k}}, \quad k \geq 1
\]
and the leading coefficient of $R_{k}(b;\, z)$ is 
\[
     \prod_{j=1}^{k} (1+ic_{j}^{(b)}) = \frac{(q^{b}; q)_{k}}{(q^{\lambda}\cos(\eta_{q}); q)_{k}}, \quad k \geq 1. 
\]

\noindent {\bf Assumption on  $b$:}\  From now on we assume that the value of $b$ be such that  $\mathcal{R}e(b) > 0$.  Then the elements of the sequence $\{d_{k+1}^{(b)}\}_{k \geq 1}$ satisfy  
\[
       d_{k+1}^{(b)} \leq  d_{k+1}^{(\lambda)},  \quad k \geq 1. 
\]
The sequence  $\{d_{k+1}^{(\lambda)}\}_{k \geq 1}$  is also the sequence of  coefficients that appear in the three term recurrence formula 
\[
    \hat{C}_{k+1}(x; q^{\lambda}\,|\,q) = x \, \hat{C}_{k}(x; q^{\lambda}\,|\,q) - d_{k+1}^{(\lambda)}\, \hat{C}_{k-1}(x; q^{\lambda}\,|\,q), \quad k \geq 1,
\]
of the  monic continuous $q$-ultraspherical polynomials $\{\hat{C}_{k}(x; q^{\lambda}\,|\,q)\}_{k \geq 0}$. These polynomials are symmetric and orthogonal on the interval $[-1,1]$. The continuous $q$-ultraspherical polynomials  were introduced by Rogers \cite{Rogers-1895} and as  recent references to these polynomials we refer to \cite{AskIsm-1982, Ismail-book}. 

Thus, the sequence $\{d_{k+1}^{(\lambda)}\}_{k \geq 1}$ is  a positive chain sequence.  This affirmation follows from a well known result regarding orthogonal polynomials defined on any finite interval of the real line that connects positive chain sequences to the values of the polynomials at extreme points of the interval of orthogonality.  Hence, by the comparison theorem for positive chain sequences \cite[p.\,97]{Chihara-book}, the sequence  $\{d_{k+1}^{(b)}\}_{k \geq 1}$ is also confirmed to be a positive chain sequence for any $b$ such that $\mathcal{R}e(b) = \lambda > 0$.  

Since $\mathcal{R}e(b)  > 0$, also from \cite{CosGodLamRan-2012}
\begin{equation} \label{Moment-Integral-Rep-q-Rn}
    \mathcal{L}^{(b)}[\zeta^{-j}] = \rho^{(b)}\, \int_{\mathbb{T}} \zeta^{-j}\frac{(q\zeta;\,q)_{\infty}\,(1/\zeta;\,q)_{\infty}} {(q^{b}\zeta;\,q)_{\infty}\,(q^{\overline{b}}/\zeta ;\,q)_{\infty}} \frac{1}{2\pi i\zeta}\,d\zeta, \quad 
\end{equation}
for $j = 0, \pm 1, \pm 2, \ldots \ $, where $   \rho^{(b)} = \frac{(q;\,q)_{\infty}\,(q^{b+\overline{b}};\,q)_{\infty}} {(q^{b};\,q)_{\infty}\,(q^{\overline{b}+1};\,q)_{\infty}}$. 
Note that when \eqref{Moment-Integral-Rep-q-Rn} holds, then we also have for $|z| < 1$,
\begin{equation} \label{Eq-qHypergeometric-Series}
     \mathcal{L}^{(b)}\Big[\frac{\zeta}{\zeta -z}\Big] = \sum_{j=0}^{\infty} \mathcal{L}^{(b)}[\zeta^{-j}]\, z^{j} = \,_2\phi_1\Big(\begin{array}{c}
                             q,\, q^{-b+1} \\[0.5ex]
                              q^{\overline{b}+1}
                             \end{array}\!\!;  \, q,\, q^{b}z\Big). 
\end{equation}

By using the three term recurrence formula \eqref{TTRR-q-Rn} together with results established in \cite{DimitarRanga-MN2013}, it was shown in \cite{AreaGodoyLamblemRanga-AMC2013} that when $\mathcal{R}e(b) > 0$ the zeros $z_{k,j}^{(b)}$, $j=1,2,\ldots,k$ of $R_{k}(b; z)$ are all {\em simple and lie on the unit circle $|z| = 1$}. Moreover, with $z_{k,j}^{(b)} = e^{i\theta_{k,j}^{(b)}}$, the  interlacing property 
\[
     0 < \theta_{k+1,1}^{(b)} < \theta_{k,1}^{(b)} < \theta_{k+1,2}^{(b)} < \cdots < \theta_{k,k}^{(b)} < \theta_{k+1,k+1}^{(b)} < 2\pi, \quad k \geq 1.
\]
also holds.

Polynomials given by a three term recurrence formula of the form 
\[
    R_{k+1}(z) = \big[(1+ic_{k+1})z+ (1-ic_{k+1})\big] R_{k}(z) - 4 d_{k+1} z\, R_{k-1}(z), \quad k \geq 1,
\]
with $R_0(z) = 1$, $R_1(z) = (1+ic_{1})z+ (1-ic_{1})$, where $\{c_k\}_{k=1}^{\infty}$ is any real sequence and $\{d_{k+1}\}_{k=1}^{\infty}$ is a positive chain sequence, have been the subject of study  in the recent publications  \cite{BracRangaSwami-2016, CastilloCostaRangaVeronese-JAT2014, CostaFelixRanga-JAT2013, DimitarRanga-MN2013, MFinkelRangaVeronese-2015}. These polynomials turn out to be  para-orthogonal polynomials related to some associated orthogonal polynomials on the unit circle. The main results of the present manuscript are obtained as applications of  results established in \cite{BracRangaSwami-2016}, \cite{CastilloCostaRangaVeronese-JAT2014} and \cite{CostaFelixRanga-JAT2013}.

To be able to apply  the results presented in  \cite{BracRangaSwami-2016}, \cite{CastilloCostaRangaVeronese-JAT2014} and \cite{CostaFelixRanga-JAT2013} directly, we now introduce the new moment functional $\mathcal{N}^{(b)}$ given by  
\begin{equation} \label{Eq-Explicit-Integral-Nb}
  \begin{array}{l}
     \dsp \frac{1}{2d_1}\mathcal{N}^{(b)}[\zeta^{-j}] =  - \frac{1 - q^{\lambda}\cos(\eta_q)}{(1-q^{\overline{b}})}\,\mathcal{L}^{(b)}[\zeta^{-j+1}] \\[2ex]
 \qquad \quad \dsp = (1 - q^{\lambda}\cos(\eta_q)) \frac{(q;\,q)_{\infty}\,(q^{b+\overline{b}};\,q)_{\infty}} {(q^{b};\,q)_{\infty}\,(q^{\overline{b}};\,q)_{\infty}} \int_{\mathbb{T}} \zeta^{-j}\frac{|(q\zeta;\,q)_{\infty}|^2} {|(q^{b}\zeta;\,q)_{\infty}|^2}\frac{1-\zeta}{2\pi i\zeta} \,d\zeta.
  \end{array}
\end{equation}
The non zero real constant $d_1$ is arbitrary. However,  in order to use results stated  in \cite{CastilloCostaRangaVeronese-JAT2014}, in Section \ref{Sec-New-OPUC2}  we will take $d_1 = (1-t)M_{1}^{(b)}$ with $0 \leq t < 1$, where $\{M_{k+1}^{(b)}\}_{k \geq 0}$ is the maximal parameter sequence of the positive chain sequence $\{d_{k+1}^{(b)}\}_{k \geq 1}$.

From \eqref{Explict-Exp-Moments-q-Rn} and \eqref{Eq-Explicit-Integral-Nb} we have
\begin{equation} \label{Eq-Explicit-moments-Nb}
    \dsp \mathcal{N}^{(b)}[\zeta^{-j}] =  \nu_{j}^{(b)} =  2d_1 \frac{1 - q^{\lambda}\cos(\eta_q)}{1-q^{b}}\, \frac{(q^{-b};\,q)_j}{(q^{\overline{b}};\,q)_j}\, q^{jb}, \quad j = 0, \pm 1, \pm 2, \ldots \ .
\end{equation}
The  moments $\nu_{j}^{(b)}$ are such 
\[
     \nu_0^{(b)} = \frac{2d_1}{1+ic_1^{(b)}} \quad \mbox{and} \quad 
     \nu_{j}^{(b)} = -\overline{\nu_{-j+1}^{(b)}}, \quad j =1, 2, 3, \ldots 
\]
Furthermore, from \eqref{Orthogonality-general-Lb} and \eqref{TTRR-q-Rn},  
\begin{equation} \label{Eq-LOrthogonality-P(b)}
    \mathcal{N}^{(b)}[\zeta^{-k+j} P_{k}(b;\zeta)] =
       \left\{ \begin{array}{rl}
          \dsp -2d_1\frac{1-q^{\lambda}\cos(\eta_q)}{1-q^{\overline{b}}}\rho_k^{(b)}, & j = -1, \\[1.5ex]
           0\ \, , & 0 \leq j \leq k-1, \\[0ex]
      \end{array} \right.
\end{equation}
and
\begin{equation} \label{Eq-LOrthogonality-R(b)}
    \mathcal{N}^{(b)}[\zeta^{-k+j} R_{k}(b;\zeta)] = 
      \left\{ \begin{array}{rl}
           -\overline{\gamma}_{k}^{(b)}, & j = -1, \\[0.5ex]
           0\ \, , & 0 \leq j \leq k-1, \\[0ex]
           \gamma_{k}^{(b)}, & j = k,
      \end{array} \right. 
\end{equation}
for $k \geq 1$, where $\gamma_{k}^{(b)} = \frac{4 d_{k+1}^{(b)}}{1+ic_{k+1}^{(b)}}\gamma_{k-1}^{(b)}$, $k \geq 1$, with $\gamma_{0}^{(b)} = \nu_{0}^{(b)}$. 

Now consider the polynomials $Q_{k}(b;z)$ defined by
\[
       Q_{k}(b;z) = \mathcal{N}^{(b)}\Big[\frac{R_{k}(b;z)- R_{k}(b;\zeta)}{z-\zeta}\Big], \quad k \geq 0.
\]
It is not difficult to show that $\{Q_k(b;z)\}_{k\geq 0}$, where  $Q_k(b;z)$ is of degree $k-1$, satisfies the three term recurrence formula
\[
   Q_{k+1}(b;\, z) = \big[(1 + i\,c_{k+1}^{(b)})z + (1 - i\,c_{k+1}^{(b)})\big]Q_{k}(b;\, z) - 4 d_{k+1}^{(b)} z Q_{k-1}(b;\, z), \quad k \geq 1,
\]
with $Q_{0}(b;\, z) = 0$ and $Q_{1}(b;\, z) = 2d_1$.  Moreover,  
\[
  \begin{array}{rl}
   \dsp   -\sum_{j=0}^{\infty} \nu_{j+1}^{(b)}z^{j} - \frac{Q_{k}(b;\, z)}{R_{k}(b;\, z)} &\dsp  = \frac{\overline{\gamma}_{k}^{(b)}}{R_{k}(b;0)}z^{k} + O(z^{k+1}), \\[3ex]
   \dsp   \sum_{j=1}^{\infty} \nu_{-j+1}^{(b)} z^{-j} - \frac{Q_{k}(b;\, z)}{R_{k}(b;\, z)} &\dsp  = \frac{\gamma_{k}^{(b)}}{\overline{R_{k}(b;0)}} \frac{1}{z^{k+1}}+  O((1/z)^{k+2})).
  \end{array}
\]

Since, 

\[
    -\sum_{j=0}^{\infty} \nu_{j+1}^{(b)}z^{j} = 2d_1\frac{1 - q^{\lambda}\cos(\eta_q)}{(1-q^{\overline{b}})} \sum_{j=0}^{\infty} \mathcal{L}^{(b)}[\zeta^{-j}]\, z^{j},
\]
from \eqref{Eq-qHypergeometric-Series} we have  for $|z| < 1$, 
\begin{equation} \label{Eq-Correspondance}
  \begin{array}{rl}
    \dsp 2d_1\frac{1 - q^{\lambda}\cos(\eta_q)}{(1-q^{\overline{b}})} \,_2\phi_1\Big(\begin{array}{c}
                             q,\, q^{-b+1} \\[0.5ex]
                              q^{\overline{b}+1}
                             \end{array}\!\!;  \, q,\, q^{b}z\Big) -  \frac{Q_{k}(b;\, z)}{R_{k}(b;\, z)} &\dsp  = \frac{\overline{\gamma}_{k}^{(b)}}{R_{k}(b;0)}z^{k} + O(z^{k+1}). 
  \end{array}
\end{equation}
Thus, we are able to state the following asymptotic results.

\begin{theo} Let $\mathcal{R}e(b) > 0$. Then 
\begin{equation} \label{Eq-Uniform-Convergence}
    \lim_ {k\to \infty} \frac{Q_{k}(b;\, z)}{R_{k}(b;\, z)} = 2d_1\frac{1 - q^{\lambda}\cos(\eta_q)}{(1-q^{\overline{b}})} \,_2\phi_1\Big(\begin{array}{c}
                             q,\, q^{-b+1} \\[0.5ex]
                              q^{\overline{b}+1}
                             \end{array}\!\!;  \, q,\, q^{b}z\Big), 
\end{equation}
uniformly on compact subsets of  $|z| < 1$.  
\begin{equation} \label{Eq-Rn-Asymp}
    \lim_{k \to \infty} R_{k}(b;z) = \frac{(q^{\overline{b}};\,q)_{\infty}}{(q^{\lambda}\cos(\eta_q);\,q)_{\infty}}\, \frac{(q^{b}z;q)_{\infty}}{(z;q)_{\infty}},
\end{equation}
uniformly on compact subsets of $|z| < 1$.  
\end{theo}

\noindent {\bf Proof}.   Since the rational functions $Q_{k}(b;\, z)/R_{k}(b;\, z)$ are analytic in $|z| < 1$, \eqref{Eq-Uniform-Convergence} is  an immediate consequence of \eqref{Eq-Correspondance}. Proof of \eqref{Eq-Rn-Asymp} follows from \eqref{Eq-Rk-ExplicitForm} by using  the Lebesque's dominated convergence theorem and $\lim_{k \to \infty} \frac{(q^{-k};q)_{j}}{(q^{-\overline{b}-k+1}; q)_{j}} = q^{(\overline{b}-1)j}$. \eProof

\setcounter{equation}{0} 
\section{Orthogonal polynomials  with respect to the measure $\hat{\mu}^{(b)}$}  \label{Sec-New-OPUC1}

We now consider the positive measure  $\hat{\mu}^{(b)}$ given by \eqref{Eq-NewMeasure-1} for $\mathcal{R}e(b) > 0$.

\begin{theo} \label{Thm-prob-constant-mu-hat}
With  $b = \lambda - i \eta$, $\eta_q = \eta \ln(q)$ and  $\lambda > 0$  let
\begin{equation} \label{Eq-tau1}
     \hat{\sigma}^{(b)} = \frac{1}{2(1 - q^{\lambda}\cos(\eta_q))} \frac{(q;\,q)_{\infty}\,(q^{b+\overline{b}};\,q)_{\infty}} {(q^{b+1};\,q)_{\infty}\,(q^{\overline{b}+1};\,q)_{\infty}}.
\end{equation}
Then the measure $\hat{\mu}^{(b)}$ given by \eqref{Eq-NewMeasure-1} is a probability measure on the unit circle. 

\end{theo}

\noindent {\bf Proof}.  To obtain information about this measure and the associated OPUC, we can apply the results obtained in \cite{BracRangaSwami-2016} with the three term recurrence formula \eqref{TTRR-q-Rn} and the associated moment functional $\mathcal{N}^{(b)}$. From the integral representation given in \eqref{Eq-Explicit-Integral-Nb} for $\mathcal{N}^{(b)}$, we observe that $\int_{\mathbb{T}} \zeta^{-j} d \hat{\mu}^{(b)}(\zeta) = const\ \mathcal{N}^{(b)}[\zeta^{-j}(1-\zeta^{-1})]$. 

Thus, using \cite[Thm.\,3.1]{BracRangaSwami-2016} we obtain $\hat{\mu}^{(b)}$ as a probability measure from  
\[
    \int_{\mathbb{T}} \zeta^{-j} d \hat{\mu}^{(b)}(\zeta) = \frac{1+(c_1^{(b)})^2}{4 d_1} \mathcal{N}^{(b)}[\zeta^{-j}(1-\zeta^{-1})], \quad j = 0, \pm1, \pm2, \ldots .
\]
Hence, first we have 
\[
    \int_{\mathbb{T}} \zeta^{-j} d \hat{\mu}^{(b)}(\zeta) = \frac{1+(c_1^{(b)})^2}{4 d_1} [\nu_{j}^{(b)} - \nu_{j+1}^{(b)}], \quad j = 0, \pm1, \pm2, \ldots ,
\]
from which, using \eqref{Eq-Explicit-moments-Nb}, one can confirm  that $\int_{\mathbb{T}}  d \hat{\mu}^{(b)}(\zeta) = 1$. 

On the other hand, using \eqref{Eq-Explicit-Integral-Nb},    
\[
 \begin{array}l
    \dsp \int_{\mathbb{T}} \zeta^{-j} d \hat{\mu}^{(b)}(\zeta) \\[2ex]
      \qquad \dsp    = \frac{1+(c_1^{(b)})^2}{2}(1 - q^{\lambda}\cos(\eta_q)) \frac{(q;\,q)_{\infty}\,(q^{b+\overline{b}};\,q)_{\infty}} {(q^{b};\,q)_{\infty}\,(q^{\overline{b}};\,q)_{\infty}} \int_{\mathbb{T}} \zeta^{-j}\frac{|(\zeta;\,q)_{\infty}|^2} {|(q^{b}\zeta;\,q)_{\infty}|^2}\frac{1}{2\pi i\zeta} \,d\zeta. 
 \end{array}
\]
Thus, from the expression for $c_1^{(b)}$ in \eqref{TTRR-q-Rn-coeffs} we arrive at the value of $\hat{\sigma}^{(b)}$. \eProof

To obtain the sequence of monic OPUC $\{\hat{\Phi}(b; z)\}_{k \geq 0}$ one needs the minimal parameter sequence of $\{d_{k+1}^{(b)}\}_{k\geq1}$. 

\begin{theo} \label{Thm-minimal-Parameter-Sequences}
With  $b = \lambda - i \eta$, $\eta_q = \eta \ln(q)$ and  $\lambda > 0$  let $\{d_{k+1}^{(b)}\}_{k\geq 1}$ be the positive chain sequence given by \eqref{TTRR-q-Rn-dn-coeffs}.
Then the minimal parameter sequence $\{\ell_{k+1}^{(b)}\}_{k\geq 0}$ of  $\{d_{k+1}^{(b)}\}_{k\geq 1}$ is such that
\[ 
   1 - \ell_{k+1}^{(b)}  = \frac{R_{k+1}(b;1)}{2R_{k}(b;1)}   =  \frac{1-q^{\overline{b}+k}}{2(1-q^{\lambda+k}\cos(\eta_{q}))}\frac{\,_2\phi_1(q^{-k-1},\, q^{b}; q^{-\overline{b} -k};q,q^{-\overline{b}+1})}{\,_2\phi_1(q^{-k},\, q^{b}; q^{-\overline{b}+1-k};q,q^{-\overline{b}+1})}, \quad k \geq 0.
\]
Moreover, $\lim_{k \to \infty} \ell_{k+1}^{(b)} = 1/2$. 

\end{theo}

\noindent {\bf Proof.} The proof of expression for $1 - \ell_{k+1}^{(b)}$ is easily obtained from rewriting the three term recurrence formula \eqref{TTRR-q-Rn} for $z=1$ in the form
\[
   d_{k+1}^{(b)} = \frac{R_{k}(b;1)}{2R_{k-1}(b;1)}\Big[1 - \frac{R_{k+1}(b;1)}{2R_{k}(b;1)}  \Big], \quad k \geq 1. 
\] 
Since $\lim_{k \to \infty} d_{k+1}^{(b)} = 1/4$, the asymptotic for $\ell_{k+1}^{(b)}$ follows from Theorem 6.4 in \cite[p.\,102]{Chihara-book}.  \eProof
 
\begin{proof}[\bf Proof of Theorem \ref{Thm-New-OPUC1}]
Let the polynomials $\hat{\Phi}_{k}(b; .)$ be given by 
\[
    \hat{\Phi}_{k}(b; z) = \frac{1}{R_{k}(b;1)\prod_{j=1}^{k+1} (1+i c_{j}^{(b)})} \frac{R_{k+1}(b;z)R_{k}(b;1)-R_{k}(b;z)R_{k+1}(b;1)}{z-1}, \quad k \geq 0.
\]
Clearly $\hat{\Phi}_{k}(b; z)$ is a monic polynomial of degree $k$. From \eqref{Eq-LOrthogonality-R(b)} and \eqref{Eq-Explicit-Integral-Nb} we can also verify that 
\[
    \int_{\T} \zeta^{-j}\hat{\Phi}_{k}(b; \zeta) d \hat{\mu}^{(b)}(\zeta) = 0, \quad j=0,1,\ldots, k-1. 
\]
Hence, $\hat{\Phi}_{k}(b; .)$ are the required monic orthogonal polynomials and with Theorem \ref{Thm-minimal-Parameter-Sequences}
\[
    \hat{\Phi}_{k}(b; z) = \frac{1}{\prod_{j=1}^{k+1} (1+i c_{j}^{(b)})} \frac{R_{k+1}(b;z) - 2(1-\ell_{k+1}^{(b)}) R_{k}(b;z)}{z-1}, \quad k \geq 0.
\]

Letting $z=0$, the associated Verblunsky coefficients are found to be   
\begin{equation} \label{Eq-Verblunski-OPUC1}
    \hat{\alpha}_{k-1}^{(b)} = -   \frac{1 - 2\ell_{k+1}^{(b)} - ic_{k+1}^{(b)}}{1 - ic_{k+1}^{(b)}}\ \prod_{j=1}^{k} \frac{1+ic_{j}^{(b)}}{1-ic_{j}^{(b)}} \, \quad k \geq 1.
\end{equation}
Hence, the first two parts of Theorem \ref{Thm-New-OPUC1} follows from \eqref{TTRR-q-Rn-coeffs}. 

To obtain the last part of Theorem \ref{Thm-New-OPUC1} we observe from \eqref{Eq-Verblunski-OPUC1} that 
\[
    1 - |\hat{\alpha}_{k-1}^{(b)}|^2 = \frac{4\ell_{k+1}^{(b)}(1-\ell_{k+1}^{(b)})}{1 + (c_{k+1}^{(b)})^2}, \quad k \geq 1. 
\]
Hence, from $(1-\ell_{j}^{(b)})\ell_{j+1}^{(b)} = d_{j+1}^{(b)}$, $j \geq 1$, where $\ell_{1}^{(b)} = 0$, we find 
\[
   \prod_{j=1}^{k} (1 - |\hat{\alpha}_{j-1}^{(b)}|^2 ) = (1-\ell_{k+1}^{(b)})\prod_{j=1}^{k}\frac{4d_{j+1}^{(b)}}{1+(c_{j+1}^{(b)})^2}, \quad k \geq 1. 
\]
Hence, the last result of Theorem \ref{Thm-New-OPUC1} follows from $(\hat{\kappa}_k^{(b)})^{-2} = \prod_{j=1}^{k}(1 - |\hat{\alpha}_{j-1}^{(b)}|^2)$. \hfill
\end{proof}

From Theorem \ref{Eq-Verblunski-OPUC1} and the reciprocal property of $R_n(b;.)$ we have  
\[
    \hat{\Phi}_{k}^{\ast}(b; z) = \frac{(q^{\lambda}\cos(\eta_{q}); q)_{k+1}}{(q^{\overline{b}};q)_{k+1}}\, \frac{R_{k+1}(b;z) - 2(1-\ell_{k+1}^{(b)})z R_{k}(b;z)}{1-z}, \quad k \geq 0.
\]
Hence, from \eqref{Eq-Rn-Asymp} 
$
   \lim_{k \to \infty} \hat{\Phi}_{k}^{\ast}(b; z) = (q^{b}z;q)_{\infty}/(z;q)_{\infty}, 
$
uniformly on compact subsets of $|z| < 1$. Thus, by considering the limit (see \cite[p.\,144]{Simon-Book-p1}) of  $1/(\hat{\kappa}_k  \hat{\Phi}_{k}^{\ast}(b; z))$ we can state   the following.

\begin{theo} The Szeg\H{o} function  associated with the probability measure $\hat{\mu}^{(b)}$ given by \eqref{Eq-NewMeasure-1} and \eqref{Eq-tau1} is 
\[
    \hat{D}(z) =  \frac{1}{|(q^{b+1};q)_{\infty}|}\sqrt{\frac{(q;q)_{\infty}(q^{2\lambda};q)_{\infty}}{2 (1-q^{\lambda}\cos(\eta_{q}))}} \, \frac{(z;q)_{\infty}}{(q^{b}z;q)_{\infty}}.
\] \\[-2ex]
\end{theo} 

\setcounter{equation}{0} 
\section{Orthogonal polynomials  with respect to the measure $\check\mu^{(b)}$}  \label{Sec-New-OPUC2}

We now look at the positive measure  $\check{\mu}^{(b)}$ given by \eqref{Eq-NewMeasure-2}.

\begin{theo} With  $b = \lambda - i \eta$, $\eta_q = \eta \ln(q)$ and  $\lambda > 0$  let 
\begin{equation} \label{Eq-tau2}
     \check{\sigma}^{(b)} = \frac{(1-q^{\overline{b}})}{\,_2\phi_1(q,\, q^{-b+1}; q^{\overline{b}+1};q,q^{b})}\, \frac{(q;\,q)_{\infty}\,(q^{b+\overline{b}};\,q)_{\infty}} {(q^{b};\,q)_{\infty}\,(q^{\overline{b}};\,q)_{\infty}}.
\end{equation}
Then the measure $\check{\mu}^{(b)}$  given by  \eqref{Eq-NewMeasure-2} is a probability measure on the unit circle. 

\end{theo}

\noindent {\bf Proof}. With  $\mathcal{R}e(b) > 0$ we clearly have 
$\,_2\phi_1(q,\, q^{-b+1}; q^{\overline{b}+1};q,q^{b})$ finite. Moreover, from  \eqref{Moment-Integral-Rep-q-Rn} and  \eqref{Eq-qHypergeometric-Series},   for $|z| < 1$, 
\begin{equation*}
  \begin{array}{l}
    \dsp \frac{1}{(1-q^{\overline{b}})} \,_2\phi_1\Big(\begin{array}{c}
                             q,\, q^{-b+1} \\[0.5ex]
                              q^{\overline{b}+1}
                             \end{array}\!\!;  \, q,\, q^{b}z\Big) 
           = \frac{1}{(1-q^{\overline{b}})}  \mathcal{L}^{(b)} \Big[\frac{\zeta}{\zeta -z}\Big], \\[3ex]
   \qquad \qquad \qquad \dsp =  \frac{(q;\,q)_{\infty}\,(q^{b+\overline{b}};\,q)_{\infty}} {(q^{b};\,q)_{\infty}\,(q^{\overline{b}};\,q)_{\infty}} \int_{\mathbb{T}} \frac{\zeta}{\zeta -z}\frac{(q\zeta;\,q)_{\infty}\,(1/\zeta;\,q)_{\infty}} {(q^{b}\zeta;\,q)_{\infty}\,(q^{\overline{b}}/\zeta ;\,q)_{\infty}}\frac{1}{2 \pi i \zeta} \,d\zeta.                            
 \end{array}
\end{equation*}
 Hence, taking the limit as $z \to 1$ from below, we have  
\[
  \begin{array}l
     \dsp \frac{\,_2\phi_1(q,\, q^{-b+1}; q^{\overline{b}+1};q,q^{b})}{(1-q^{\overline{b}})}  
       = \frac{(q;\,q)_{\infty}\,(q^{b+\overline{b}};\,q)_{\infty}} {(q^{b};\,q)_{\infty}\,(q^{\overline{b}};\,q)_{\infty}} \int_{\mathbb{T}} \frac{(q\zeta;\,q)_{\infty}\,(q/\zeta;\,q)_{\infty}} {(q^{b}\zeta;\,q)_{\infty}\,(q^{\overline{b}}/\zeta ;\,q)_{\infty}}\frac{1}{2\pi i\zeta} \,d\zeta \, > 0.
   \end{array}
\]
The above result can be justified by Abel's continuity theorem. Thus, with $\check{\sigma}^{(b)}$ given as in the theorem, we have 
\[
   \check{\sigma}^{(b)} \int_{\mathbb{T}} \frac{(q\zeta;\,q)_{\infty}\,(q/\zeta;\,q)_{\infty}} {(q^{b}\zeta;\,q)_{\infty}\,(q^{\overline{b}}/\zeta ;\,q)_{\infty}}\frac{1}{2\pi i\zeta} \,d\zeta = 1.  \\[-3ex]
\]
\eProof

From now on we will assume that $\check{\mu}^{(b)}$ is a probability measure. Hence,  from \eqref{Eq-Explicit-Integral-Nb} it is not difficult to verify that 
\begin{equation} \label{Eq-Explicit-Relation-Lb-muCheck}
   \begin{array}{l}
      \dsp \frac{-1}{2d_1} \frac{(1-q^{\overline{b}})}{1 - q^{\lambda}\cos(\eta_q)} \frac{\mathcal{N}^{(b)}[\zeta^{-j}]}{\,_2\phi_1(q,\, q^{-b+1}; q^{\overline{b}+1};q,q^{b})} 
     = \frac{\mathcal{L}^{(b)} [\zeta^{-j+1}]}{\,_2\phi_1(q,\, q^{-b+1}; q^{\overline{b}+1};q,q^{b})} \\[3ex]
      \qquad \qquad \qquad \qquad \qquad \quad \dsp=  \int_{\mathbb{T}} (1-\zeta^{-1})\zeta^{-j+1} d \check{\mu}^{(b)}(\zeta), \quad j =0, \pm1, \pm2, \ldots \ .
   \end{array}
\end{equation}

Now let us consider the monic polynomials 
\[
    A_{k}(z) = \frac{\check{\Phi}_{k}(b; z) - \check{\tau}_k^{(b)}\,\check{\Phi}_{k}^{\ast}(b;z)}{z-1}, \quad k \geq 1,
\] 
where
\[
   \check{\tau}_k^{(b)} = \frac{\check{\Phi}_{k}(b; 1)}{\check{\Phi}_{k}^{\ast}(b; 1)}, \quad k \geq 0.
\]

The polynomial $A_{k}(z)$ is a constant multiple of $\check{K}_{k}(b; z, 1)$, where  
\[
      \check{K}_{k}(b; z, w) = \sum_{j=0}^{k} \overline{\check{\phi}_{j}(b; w)}\, \check{\phi}_{j}(b; z),  \quad k \geq 0, 
\] 
are  the associated Christoffel-Darboux (or CD)  kernels. From the orthogonality of \linebreak $\{\check{\Phi}_{k}(b; z)\}_{k \geq 0}$ it is easy to see that 
\[
    \int_{\mathbb{T}} \zeta^{-k+j} A_{k}(\zeta) (1-\zeta) d\check{\mu}^{(b)}(\zeta) = 0, \quad 0 \leq j \leq k-1.
\]
Hence, from \eqref{Eq-Explicit-Relation-Lb-muCheck} we also have 
\[
     \mathcal{N}^{(b)}[\zeta^{-k+j} A_{k}(\zeta)] = 0, \quad 0 \leq j \leq k-1.
\]
Thus, by comparing the determinant representation for the monic polynomials  $A_{k}(z)$ obtained  from the above orthogonality conditions with the determinant representation for the monic polynomials $P_{k}(b;z)$ obtained  from \eqref{Eq-LOrthogonality-P(b)},  it follows that 
\[
   A_{k}(z) = P_{k}(b;z), \quad k \geq 0.
\]  
Moreover, comparing the three term recurrence formula in  \cite[Thm.\,2.2]{CostaFelixRanga-JAT2013} for the polynomials 
\[
    \frac{\prod_{j=0}^{k-1} [1-\check{\tau}_{j}^{(b)} \check{\alpha}_j^{(b)}]}{\prod_{j=0}^{k-1} [1-\Re(\check{\tau}_{j}^{(b)} \check{\alpha}_j^{(b)})]}\, A_{k}(z), \quad k \geq 1,
\]
with the three term recurrence formula \eqref{TTRR-q-Rn}, we have 
\[
    \frac{-\Im[\check{\tau}_{k-1}^{(b)}\check{\alpha}_{k-1}^{(b)}]}{1-\Re[\check{\tau}_{k-1}^{(b)}\check{\alpha}_{k-1}^{(b)}]} = c_{k}^{(b)}, \quad k \geq 1
\]
and
\[
     \big(1-g_{k}^{(b)}\big)\,g_{k+1}^{(b)} = d_{k+1}^{(b)}, \quad k \geq 1,
\]
with 
\[
    g_{k}^{(b)} = \frac{1}{2} \frac{|1- \check{\tau}_{k-1}^{(b)} \check{\alpha}_{k-1}^{(b)}|^2}{1-\Re[\check{\tau}_{k-1}^{(b)} \check{\alpha}_{k-1}^{(b)}]}.
\]
The sequence $\{g_{k+1}^{(b)}\}_{k\geq 0}$ is also shown in \cite{CostaFelixRanga-JAT2013} to be a parameter sequence for the positive chain sequence $\{d_{k+1}^{(b)}\}_{k\geq 1}$.  Moreover, since the measure $\check{\mu}^{(b)}$ does not have a pure mass point (or pure point) at $z=1$,  the sequence  $\{g_{k+1}^{(b)}\}_{k\geq 0}$ is the maximal parameter sequence of $\{d_{k+1}^{(b)}\}_{k\geq 1}$. In what follows, we denote $g_{k}^{(b)} = M_{k}^{(b)}$, for $k \geq 1$. Hence, 
\[ 
    M_{1}^{(b)} =\frac{1}{2} \frac{|1- \check{\alpha}_{0}^{(b)}|^2}{1-\mathcal{R}e[\check{\alpha}_{0}^{(b)}]} .
\]
To find an explicit expression for $M_{1}^{(b)}$ we need to evaluate $\check{\alpha}_{0}^{(b)}$.
Observe that $\check{\alpha}_{0}^{(b)} = \check{\mu}_1^{(b)} = \int_{\mathbb{T}} \zeta^{-1} d \check{\mu}^{(b)}(\zeta)$.  Hence, from $\zeta^{-1} = 1 - (1-\zeta^{-1})$, we obtain 
\[ 
  \begin{array}{ll}
    \check{\alpha}_{0}^{(b)} &\dsp = 1 - \int_{\mathbb{T}} (1-\zeta^{-1}) d \check{\mu}^{(b)}(\zeta)\\[2ex]
     & = 1 - \mathcal{L}^{(b)}[1]/\,_2\phi_1(q,\, q^{-b+1};     q^{\overline{b}+1};q,q^{b}).   
  \end{array}
\]
Since  $\mathcal{L}^{(b)}[1] = 1$ we then have
\[
     \check{\alpha}_{0}^{(b)} = 1 - [\,_2\phi_1(q,\, q^{-b+1}; q^{\overline{b}+1};q,q^{b})]^{-1}. 
\]
From this and by  observing from \eqref{Eq-tau2} that $(1- q^{\overline{b}})^{-1}\,_2\phi_1(q,\, q^{-b+1}; q^{\overline{b}+1};q,q^{b})$ is also positive, we obtain 
\begin{equation} \label{Eq-Initial-Max-Param-Seq}
    M_{1}^{(b)} = \frac{1}{2} \frac{1- q^{\overline{b}}}{1- q^{\lambda}\cos(\eta_q)} \frac{1}{\,_2\phi_1(q,\, q^{-b+1}; q^{\overline{b}+1};q,q^{b})}. 
\end{equation}

\begin{theo} \label{Thm-maximal-Parameter-Sequences}
With  $b = \lambda - i \eta$, $\eta_q = \eta \ln(q)$ and  $\lambda > 0$   let $\{d_{k+1}^{(b)}\}_{k\geq 1}$ be the positive chain sequence given by \eqref{TTRR-q-Rn-dn-coeffs}. 
Then if $\{M_{k+1}^{(b)}\}_{k\geq 0}$ is the maximal parameter sequences  of  $\{d_{k+1}^{(b)}\}_{k \geq 1}$ then  
\[
     M_{k+1}^{(b)} =  \frac{1}{2}\frac{1-q^{\overline{b}+k}}{1-q^{\lambda+k}\cos(\eta_{q})} \frac{\,_2\phi_1(q^{k},\, q^{-b+1}; q^{\overline{b}+k};q,q^{b})}{\,_2\phi_1(q^{k+1},\, q^{-b+1}; q^{\overline{b}+k+1};q,q^{b})}, \quad k \geq 0.
\]
Moreover, $\lim_{k \to \infty} M_{k+1}^{(b)} = 1/2$.

\end{theo}

\noindent{\bf Proof}. The value of $M_{1}^{(b)} $ is confirmed from \eqref{Eq-Initial-Max-Param-Seq}.  To verify the value of $M_{k+1}^{(b)}$ for $k \geq 1$, we consider the   relation  
\begin{equation*} \label{Eq-Contiguous-Relation-2}
    \frac{\mathfrak{C}_{k}^{(b)}}{f_{k-1}^{(b)}(z)} = z +\mathfrak{C}_{k}^{(b)} - \dsp \mathfrak{D}_{k+1}^{(b)}  z \frac{f_{k}^{(b)}(z)}{\mathfrak{C}_{k+1}^{(b)}}, \quad k \geq 1,
\end{equation*}
where $\mathfrak{C}_k^{(b)}$ and $\mathfrak{D}_{k+1}^{(b)}$ are as in \eqref{TTRR-q-Pn} and 
\[ 
    f_{k}^{(b)}(z) =  \frac{\,_2\phi_1(q^{k+1}, q^{-b+1};\, q^{\overline{b}+k+1};\,q,q^{b}z)}{\,_2\phi_1(q^{k}, q^{-b+1};\, q^{\overline{b}+k};\,q,q^{b}z)},  \quad k \geq 1.   
\]
The above relation follows from contiguous formulas  for basic hypergeometric polynomials (see \cite[p. 22]{GasRah-book}) obtained by Heine. 

Thus, we have 
\[
   \begin{array}l
   \dsp \Big(1 - \frac{\mathfrak{C}_{k}^{(b)}}{(1+\mathfrak{C}_{k}^{(b)}) f_{k-1}^{(b)}(1)}\Big) \frac{\mathfrak{C}_{k+1}^{(b)}}{(1+\mathfrak{C}_{k+1}^{(b)}) f_{k}^{(b)}(1)}  = \frac{\mathfrak{D}_{k+1}^{(b)}}{(1+\mathfrak{C}_{k}^{(b)})(1+\mathfrak{C}_{k+1}^{(b)})} , \quad k \geq 1.
   \end{array}
\]
Hence, observing that 
\[
    \frac{\mathfrak{D}_{k+1}^{(b)}}{(1+\mathfrak{C}_{k}^{(b)})(1+\mathfrak{C}_{k+1}^{(b)})} = d_{k+1}^{(b)}, \quad k \geq 1, 
\]
and 
\[
    \frac{\mathfrak{C}_{k+1}^{(b)}}{(1+\mathfrak{C}_{k+1}^{(b)}) f_{k}^{(b)}(1)} = \frac{1}{2}\frac{1-q^{\overline{b}+k}}{1-q^{\lambda+k}\cos(\eta_{q})} \frac{\,_2\phi_1(q^{k}, q^{-b+1};\, q^{\overline{b}+k};\, q,q^{b})}{\,_2\phi_1(q^{k+1}, q^{-b+1};\,  q^{\overline{b}+k+1};\, q,q^{b})},
\]
for $k \geq 0$, we obtain $(1-M_{k}^{(b)})M_{k+1}^{(b)} = d_{k+1}^{(b)}$, $k \geq 1$, which gives the expression for $M_{k+1}^{(b)}$. Now to arrive at the limit for $M_{k+1}^{(b)}$ we again use Theorem 6.4 in \cite[p.\,102]{Chihara-book}. \eProof

From results given in  \cite{CastilloCostaRangaVeronese-JAT2014}, the sequence $\{Q_{k}(b;\, 1)/[2d_1R_{k}(b;\, 1)]\}$ is a positive and increasing sequence, and that   
\[
     \lim_{k\to \infty} \frac{Q_{k}(b;\, 1)}{R_{k}(b;\, 1)} = \frac{2d_1}{2M_1^{(b)}},
\]
where $\{M_{k+1}^{(b)}\}_{k\geq 0}$ is the maximal parameter sequence of the positive chain sequence $\{d_{k+1}^{(b)}\}_{k\geq 1}$. Thus, from \eqref{Eq-Uniform-Convergence}, by continuity %
\[
     \frac{1 - q^{\lambda}\cos(\eta_q)}{(1-q^{\overline{b}})} \,_2\phi_1\Big(\begin{array}{c}
                             q,\, q^{-b+1} \\[0.5ex]
                              q^{\overline{b}+1}
                             \end{array}\!\!;  \, q,\, q^{b}\Big) =  \frac{1}{2M_1^{(b)}}, 
\]                             
which again confirms \eqref{Eq-Initial-Max-Param-Seq}. 

With the choice $d_1 = M_1^{(b)}$, we also have from \eqref{Eq-Explicit-Relation-Lb-muCheck}
\[
   \mathcal{N}^{(b)}[\zeta^{j}] = \int_{\mathbb{T}} (1-\zeta)\,\zeta^{-j} d \check{\mu}^{(b)}(\zeta), \quad j =0, \pm1, \pm2, \ldots .
\] 

Now,  with the three term recurrence formula \eqref{TTRR-q-Rn} and the associated moment functional $\mathcal{N}^{(b)}$, application of the results given in \cite{CastilloCostaRangaVeronese-JAT2014} and \cite{CostaFelixRanga-JAT2013} gives the following.

For $0 \leq t < 1$, let $\check{\mu}^{(b)}(t; .)$ be the probability measure given by $(1-t)\check{\mu}^{(b)}(.)+ t \delta_1$. Let $\{m_k^{(b,t)}\}_{k\geq 0}$ be the minimal parameter sequence of the positive chain sequence $\{d_1^{(b,t)}, d_2^{(b)}, d_3^{(b)}, d_4^{(b)}, \ldots\}$, where $d_1^{(b,t)} = (1-t)M_1^{(b)}$, then the monic OPUC $\check{\Phi}_k(b,t; z)$ associated with $\check{\mu}^{(b)}(t; .)$ are given by
\[
     \check{\Phi}_k(b,t; z) = \frac{1}{\prod_{j=1}^{k}(1+ic_j^{(b)})} [R_{k}(b;z) - 2(1-m_k^{(b,t)})R_{k-1}(b;z)], \quad k \geq 1.
\]
Observe that $m_k^{(b,0)} = M_{k}^{(b)}$, $k \geq 1$. \\[-2ex]

\begin{proof}[\bf Proof of Theorem \ref{Thm-New-OPUC2}]
To obtain the first part of Theorem \ref{Thm-New-OPUC2} we simply set $t=0$ in the above formula  and then use \eqref{TTRR-q-Rn-coeffs}. The result for the associated Verblunsky coefficients is obtained with the substitution $z =0$. To obtain the expression for $M_{1}^{(b)}$ in the theorem we use 
\[
   2 M_{1}^{(b)} = \frac{1 - q^{\overline{b}}}{1 - q^{\lambda}\cos(\eta_{q})}(1 - \check{\alpha}_{0}^{(b)})= \frac{1 - q^{b}}{1 - q^{\lambda}\cos(\eta_{q})}(1 - \overline{\check{\alpha}_{0}^{(b)}}), 
\]
which follows form the result for the Verblunsky coefficients in the Theorem. 

Now to obtain the last part of Theorem \ref{Thm-New-OPUC2},  observe that 
\[
    1 - |\check{\alpha}_{k-1}^{(b)}|^2 = \frac{4M_{k}^{(b)}(1-M_{k}^{(b)})}{1 + (c_{k}^{(b)})^2}, \quad k \geq 1. 
\]
Hence, the result follows from $(\check{\kappa}_k^{(b)})^{-2} = \prod_{j=1}^{k}(1 - |\check{\alpha}_{j-1}^{(b)}|^2)$ and $(1-M_{k}^{(b)})M_{k+1}^{(b)} = d_{k+1}^{(b)}$ for $k \geq 1$. 
\end{proof}

Now from 
\[
    \check{\Phi}_{k}^{\ast}(b;z) =  \frac{(q^{\lambda}\cos(\eta_{q}); q)_{k}}{(q^{\overline{b}};q)_{k}}\, [R_{k}(b;z) - 2(1-M_k^{(b)})zR_{k-1}(b;z)], \quad k \geq 1,
\]
we find
$
    \lim_{k \to \infty} \check{\Phi}_{k}^{\ast}(b;z) =  (q^{b}z;q)_{\infty}/(qz;q)_{\infty}, 
$
uniformly on compact subsets of $|z| < 1$. Thus, by considering the limit (see \cite[p.\,144]{Simon-Book-p1})  of  $1/(\check{\kappa}_k  \check{\Phi}_{k}^{\ast}(b; z)$ we can state  the following.

\begin{theo} If $\check{D}(z)$ is the Szeg\H{o} function  associated with the probability measure $\check{\mu}^{(b)}$ given by \eqref{Eq-NewMeasure-2} and \eqref{Eq-tau2} then
\[
    \check{D}(z) = \frac{\sqrt{2 (q;q)_{\infty}(q^{2\lambda};q)_{\infty}(1-q^{\lambda}\cos(\eta_{q}))M_1^{(b)}}}{|(q^{b};q)_{\infty}|} \frac{(qz;q)_{\infty}}{(q^{b}z;q)_{\infty}}.
\] \\[-2ex]

\end{theo}


\end{document}